\documentclass[10pt]{amsart}
\usepackage{amsmath}
\usepackage{amsthm}
\usepackage{amscd}
\usepackage{amsfonts}
\usepackage{amssymb}
\usepackage{latexsym}
\newtheorem{thm}{Theorem}
\newtheorem{prop}[thm]{Proposition}
\newtheorem{cor}[thm]{Corollary}
\newtheorem{lem}[thm]{Lemma}

\newenvironment{prf}{\begin{proof}[Proof]}{\end{proof}}

\usepackage{amsfonts}
\newcommand{\field}[1]{\mathbb{#1}}
\newcommand{\Q}{\field{Q}}

\newcommand{\N}{\field{N}}
\newcommand{\Z}{\field{Z}}

\newcommand{\F}{\field{F}}

\newcommand{\dc}{{\bf C}_{\infty}}
\newcommand{\cp}{{\bf C}_{\frak p}}

\newcommand{\Gal}{\text{Gal}}

\newcommand{\T}{\mathcal T}

\newcommand{\Hom}{\text{Hom}}

\newcommand{\n}{\mathcal N}

\newcommand{\pr}{\frak p}
\newcommand{\m}{\frak m}

\newcommand{\mm}{\frak M}
\newcommand{\ol}{\mathcal O}
\newcommand{\liminv}{\displaystyle \lim_{\leftarrow}}
\newcommand{\limdir}{\displaystyle \lim_{\rightarrow}}
\newcommand{\dlog}{{\rm dlog\,}}
\newcommand{\id}{{\rm id}}
\newcommand{\psqrt}{\;{}_{\eta}\!\!\sqrt}

\begin{document}

\title[Coleman's power series]{Coleman's power series and Wiles' reciprocity for rank 1 Drinfeld modules}

\author{Francesc Bars, Ignazio Longhi}

\thanks{Both authors have been supported by TMR Arithmetic Algebraic Geometry. The first author is also supported by MTM2006-11391; the second one by a scholarship of Universit\`a di Milano.}

\date{\today}

\maketitle

\begin{center}
\begin{small}
\begin{abstract}
We introduce the formalism of Coleman's power series for rank 1 Drinfeld modules and apply it to formulate and prove the analogue of Wiles' explicit reciprocity law in this setting.
\end{abstract}
\end{small}
\end{center}

\section{Introduction}

Andrew Wiles discovered an explicit reciprocity law for local fields \cite{wiles}, generalizing earlier work of Artin-Hasse \cite{artinhasse} and Iwasawa \cite{iwasawa}.  Since the publication of \cite{wiles}, different proofs of this reciprocity law were found: see the expositions in \cite[chapter 9]{lang}, \cite[I, \S4]{desh} (where the main tool are Coleman's power series) and \cite[\S3.3]{scholl} (a cohomological approach, inspired by Kato's formulation in \cite[\S1]{kato}).

Wiles' explicit reciprocity law has its foundation in the theory of Lubin-Tate formal groups: torsion points of the formal group generate a ``cyclotomic" tower of local fields and the action of local norm symbols on the torsion is expressed by an analytic formula.

A rank 1 Drinfeld module of generic characteristic can be seen as originating a special instance of Lubin-Tate formal group: therefore there should be an analogue of Wiles' reciprocity law in this setting. For the Carlitz module, this was proven by Angl\`es: his paper \cite{angles} follows the approach of \cite[chapter 9]{lang}.\\

In the present paper we introduce the formalism of Coleman's power
series in the additive setting of (formal) rank 1 Drinfeld
modules. As a first application we obtain our main result, the
explicit reciprocity law for Drinfeld modules over any global
field $F$ (Theorem \ref{reclaw}): i.e., we work with any
sign-normalized rank 1 Drinfeld module, with no restriction on the
class number. Following \cite[I, \S4]{desh} we formulate the law
directly in its limit form, as the equality of two pairings on
systems of groups.

In the classical situation over $\Q_p$, one can exploit the Coleman isomorphism (corresponding to Theorem \ref{teo8} in our setting) to construct $p$-adic $L$-functions: e.g., the Kubota-Leopoldt zeta function comes from the system of cyclotomic units. We hope to be able to recover this aspect of the theory in a future paper.\\

Let us finally give a rough sketch of the contents of this paper. In \S\ref{setting} we introduce the basic notation and properties for formal Drinfeld modules and the ``cyclotomic" tower we are working with. In \S\ref{coleman} we rewrite Coleman's formalism in the setting of Drinfeld modules: our construction is quite detailed, in the hope it can be a good introduction to the subject (the differences with characteristic 0 are irrelevant). In \S\ref{law} we introduce the two pairings and prove their equality.

\section{Setting}\label{setting}

Let $F$ be a global function field, with field of constants $\F_q$, $q$ a power of $p$. Given a place $v$ of $F$, $q_v:=q^{\deg(v)}$ will indicate the cardinality of its residue field $\F_v$. Besides, we denote by $\tau$ the operator $x\mapsto x^q$ and, for an $\F_p$-algebra $R$, by $R\{\tau\}$ the ring of skew polynomials with coefficients in $R$: multiplication in $R\{\tau\}$ is given by composition.

\subsection{Review of rank 1 Drinfeld modules} We briefly recall Hayes' theory of explicit construction of class fields by means of rank 1 Drinfeld modules; our main references will be \cite{hayes} and \cite[chapter 7]{goss}. We fix a place $\infty$ of $F$ and let $A\subset F$ be the ring of functions regular away from $\infty$.

As in \cite{hayes} and \cite{goss} we fix a sign-function $sgn:F_{\infty}\rightarrow\F_{\infty}$: then the basic extensions of $F$ are $H$ and $H^+$, the Hilbert class field and the normalizing field (see e.g.\! \cite[\S14-15]{hayes} or \cite[\S7.1 and \S7.4]{goss}), with $B$ and $B^+$ the integral closure of $A$ in $H$ and $H^+$ respectively. (It might be worth to recall that $H$ depends only on the choice of $\infty$, while $H^+$ is determined by both $\infty$ and $sgn$.) Let $\Phi$ be a $sgn$-normalized rank 1 Drinfeld $A$-module: i.e., $\Phi$ is a ring homomorphism $A\rightarrow B^+\{\tau\}$, $a\mapsto\Phi_a$, such that the constant coefficient of $\Phi_a$ is $a$ and the leading coefficient map is a $Gal(\F_{\infty}/\F_q)$-twist of $sgn$. We fix such a $\Phi$.

As usual, if $I$ is an ideal of $A$, $\Phi[I]$ denotes the $I$-torsion of $\Phi$ (i.e., the common zeroes of all $\Phi_a$, $a\in I$) and $\Phi_I$ is the unique monic generator of the left ideal of $H^+\{\tau\}$ generated by $\Phi_a$, $a\in I$. One sees immediately that
$$\Phi_I(x)=\prod_{u\in\Phi[I]}(x-u)$$
and since elements in $\Phi[I]$ are all integral above $B^+$, it follows that $\Phi_I\in B^+\{\tau\}$. By \cite[Proposition 11.4]{hayes}, in case $I=\pr^n$, $\pr$ a prime, all irreducible factors (over $B^+$) of the polynomial  $\Phi_I(x)$ are Eisenstein at $\pr$.

To conclude, we recall that the extension $H^+(\Phi[\pr^n])/F$ is abelian and
$$Gal(H^+(\Phi[\pr^n])/H^+)\simeq(A/\pr^n)^{\ast}$$
(\cite[\S7.5]{goss}): the isomorphism is given by the $A$-action on $\Phi[\pr^n]$.\\

\noindent{\em Caveat.} The notation $\Phi_a$ (or $\Phi_I$) will be used to denote both the operator $\Phi_a\in B^+\{\tau\}$ and the polynomial $\Phi_a(x)\in B^+[x]$; the context should make clear which one we mean.\\

\subsection{Our local setting} From now on, we fix a prime ideal $\pr$ in $A$: since $\Phi$ is $sgn$-normalized, it has good reduction in $\pr$. Let $F_{\pr}$ and $A_{\frak p}$ be the completions at $\pr$. We also fix ${\bf C}_{\pr}$, completion of an algebraic closure of $F_{\pr}$, and choose an embedding $H^+\hookrightarrow{\bf C}_{\pr}$; let $K$ be the topological closure of $H^+$ in ${\bf C}_{\pr}$ and $\ol=\ol_K$ its ring of integers, with maximal ideal $\m=\m_K$. Let $\mm:=B(0,1)=\{z\in\cp:|z|<1\}$. In what follows, all extensions of $F$ are assumed to be contained in ${\bf C}_{\pr}$. The valuation $v$ on $\cp$ is normalized so that $v(F_{\pr}^{\ast})=\Z$.

The extension $H^+/F$ is unramified outside of $\infty$: in particular, it follows that $K/F_{\pr}$ is unramified.\\

Let $\mathcal I$ be the group of fractional ideals of $A$ and $\mathcal P^+$ the subgroup of the positively generated principal ones. By Hayes' theory we know that $\mathcal I/\mathcal P^+\simeq Gal(H^+/F)$ (\cite[Theorem 14.7]{hayes}); this isomorphism is given by the action of ideals on the set of $sgn$-normalized rank 1 Drinfeld modules. In particular the image of $\pr$ in $\mathcal I/\mathcal P^+$ corresponds to the Frobenius at $\pr$: its order is $f:=[K:F_{\pr}]$, hence we have $\pr^f=(\eta)$, with $\eta\in A$ a positive element (i.e., $sgn(\eta)=1$). Then, by definition of $sgn$-normalized, $\Phi_{\eta}$ is monic and therefore $\Phi_{\eta}=\Phi_{\pr^f}$. Notice that $\eta$ is uniquely defined once we have fixed $sgn$.

Finally, we remark that, since all factors of $\Phi_{\pr^n}$ are Eisenstein, $\Phi[\pr^n]\subset\mm$ for all $n$. In particular, all coefficients (but the leading one) of $\Phi_{\eta}$ are in $\m$.

\subsection{The formal module} Consider the ring of skew power series $\ol\{\{\tau\}\}$: it is a local ring, complete in the topology induced by its maximal ideal $\m+\ol\{\{\tau\}\}\tau$. As observed by Rosen (\cite[p.\,247]{rosen}), the homomorphism $\Phi:A\rightarrow B^+\{\tau\}\subset\ol\{\{\tau\}\}$ can be extended to the localization of $A$ at $\pr$ and then to its completion: we get a formal Drinfeld module $\Phi:A_{\pr}\rightarrow\ol\{\{\tau\}\}$.

\begin{prop}[Rosen] \label{lambda} There exists a unique $\lambda\in K\{\{\tau\}\}$ of the form $\lambda=1+...$ and such that $a\lambda=\lambda\Phi_a$ for all $a\in A_{\pr}$. Besides, $\lambda$ converges on $\mm$ and it enjoys the following properties:\begin{enumerate}
\item $\lambda=\sum c_i\tau^i$ with $v(c_i)\geq -i$;
\item if $v(x)>(q-1)^{-1}$, then $v(\lambda(x))=v(x)$.
\end{enumerate}
\end{prop}

For the proof the reader is referred to \cite[Proposition 2.1 and Proposition 2.3]{rosen}.

\subsection{The ``cyclotomic'' tower} \label{tower} Define the tower of field extensions $K_n/K$ by $K_n:=K(\Phi[\pr^{fn}])$, $n\geq 1$, with ring of integers $\ol_n$ and maximal ideal $\frak m_n$. (This choice of indexing seemed to us more convenient even if slightly unusual - commonly, one calls the tamely ramified extension $K_0$.) Let $Tr_m^n:K_n\rightarrow K_m$ and $N_m^n:K_n^*\rightarrow K_m^*$ denote respectively trace and norm.\\
We also put $K_{\infty}:=\cup_{n=1}^{\infty}K_n\,$.\\

Consider the Tate module $T_{\pr}\Phi:=\liminv \Phi[\pr^{fn}]$ (the limit is taken with respect to $x\mapsto\Phi_{\eta}(x)$).\footnote{More canonically, $T_{\pr}\Phi$ is usually defined as $\Hom_{A_{\pr}}(F_{\pr}/A_{\pr}, \limdir \Phi[\pr^n])$ (\cite[Definition 4.10.9]{goss} and following remarks), but the two are isomorphic and our definition suits better our purpose.} The ring $A_{\pr}$ acts on $T_{\pr}\Phi$ via $\Phi$: i.e., $a\cdot u:=\Phi_a(u)$. Since $\Phi$ has rank 1, $T_{\pr}\Phi$ is a free $A_{\pr}$-module of rank 1.

Let $\omega=\{\omega_n\}_{n\geq 1}$ be a generator: $T_{\pr}\Phi=A_{\pr}\cdot\omega$. This means that the sequence $\{\omega_n\}$ satisfies
$$\Phi_{\eta}^n(\omega_n)=0\neq\Phi_{\eta}^{n-1}(\omega_n)\text{ and }\Phi_{\eta}(\omega_{n+1})=\omega_n.$$
(Here and in the following, when we write  $\Phi_a^n$ the power is always taken in $\ol\{\tau\}$ - that is, with respect to composition.)

By definition $K_n=K(\omega_n)$. Being a root of an Eisenstein polynomial, $\omega_n$ is a uniformizer for the field $K_n$: it follows that the extensions $K_n/K$ are totally ramified, $Gal(K_n/K)\simeq(A/\pr^{fn})^{\ast}$ and $\ol_n=\ol[[\omega_n]]=\ol[\omega_n]\,$.\\

The Galois action on $T_{\pr}\Phi$ is via the ``Carlitz-Hayes" character $\chi:G_K\rightarrow A_{\pr}^{\ast}$, defined by $\sigma\omega=\chi(\sigma)\cdot\omega$ for $\sigma\in G_K:=Gal(K_{sep}/K)$: that is, $\chi(\sigma)$ is the unique element in $A_{\pr}^{\ast}$ such that $\Phi_{\chi(\sigma)}(\omega_n)=\sigma\omega_n$ for all $n$.

\begin{lem} \label{lema2} The elements $\omega_n$ form a compatible system under the norm maps:
$$N_n^{n+1}(\omega_{n+1})=\omega_n.$$
\end{lem}

\begin{prf} Either $[K_{n+1}:K_n]$ is odd or the characteristic is 2: in both cases, it follows that $-N_n^{n+1}(\omega_{n+1})$ is the constant term of the minimal polynomial of $\omega_{n+1}$ on $K_n$. It is immediate to see that the latter is $\Phi_{\eta}(x)-\omega_n$ (it is Eisenstein).
\end{prf}

We also remark that the $Gal(K_{n+i}/K_n)$-orbit of $\omega_{n+i}$
is exactly $\omega_{n+i}+\Phi[\pr^{if}]$.

\begin{lem} \label{lema3} Let $\frak D_n:=\frak D_{K_n/F_{\pr}}$ be the different of $K_n/F_{\pr}$: then $\frak D_n$ is generated by an element of valuation
$$nf-\frac{1}{q_{\pr}-1}.$$
Moreover, $\omega_n$ has valuation
$$v(\omega_n)=[K_n:K]^{-1}=q_{\pr}^{1-nf}(q_{\pr}-1)^{-1}.$$
\end{lem}

\begin{prf} The last assertion is obvious from the already remarked fact that $\omega_n$ is a uniformizer in a totally ramified extension. As for the first, let $\psi_n$ be the irreducible polynomial of $\omega_n$ over $H^+$, $n\geq 1$: then $\frak D_n=\frak D_{K_n/K}$ because $K/F_{\pr}$ is unramified and since $\frak D_{K_n/K}=(\psi_n'(\omega_n))$ we just need to compute this derivative.\\
From \cite[Proposition 11.4]{hayes} we get the equality in $B^+[x]$
$$\psi_n(x)\Phi_{\pr^{fn-1}}(x)=\Phi_{\pr^{fn}}(x)=\Phi_{\eta^n}(x)\,.$$
Differentiating and evaluating in $\omega_n$ we get
$$v(\psi_n'(\omega_n))=nf-v(\Phi_{\pr^{fn-1}}(\omega_n))=nf-(q_{\pr}-1)^{-1}.$$
(Observe that $\Phi_{\pr^{fn-1}}(\omega_n)$ has valuation $(q_{\pr}-1)^{-1}$ because $\Phi_{\pr^{fn-1}}(x)$ is a monic polynomial of degree $q_{\pr}^{nf-1}$ all whose coefficients but the leading one are in $\m_K$.)
\end{prf}

\begin{cor} \label{trace} Let $Tr_n:K_n\rightarrow F_{\pr}$ be the trace map. Then
$$v(Tr_n(x))\geq \lfloor v(x)+nf-(q_{\pr}-1)^{-1}\rfloor$$
where $\lfloor r\rfloor$ denotes the largest integer $\leq r$. Similarly, for $m\geq 1$,
$$v\left(Tr^n_{m}(x)\right)> v(x)+(n-m)f-v(\omega_m)\,.$$
\end{cor}

\begin{prf} Let $k=\lfloor v(x)+nf-(q_{\pr}-1)^{-1}\rfloor\,$: then $x\ol_n\subseteq\pr^k\frak D_n^{-1}$ by Lemma \ref{lema3} and this means that $Tr_n(x\ol_n)\subseteq\pr^k A_{\pr}$, by a basic property of the different (see e.g. \cite[III, \S3, Proposition 7]{serre}).\\
In the same way, using the fact (obvious from Lemma \ref{lema3})
that the generator of $\frak D_{K_n/K_m}$ has valuation $(n-m)f$,
one gets
$$v\left(Tr^n_{m}(x)\right)\geq\big\lfloor\frac{v(x)+(n-m)f}{v(\omega_m)}\big\rfloor v(\omega_m)\,;$$
the second statement follows.
\end{prf}

\subsection{Local class field theory} \label{lcft} The Carlitz-Hayes character induces an isomorphism of topological groups
$$\chi^{-1}:A_{\pr}^*\rightarrow Gal(K_{\infty}/K)$$
(recall that we put $K_{\infty}=\cup_{n=1}^{\infty}K_n$). This should be compared with the local norm symbol map
$$(\,\cdot\,,K_{\infty}/F_{\pr})\colon F_{\pr}^*\rightarrow Gal(K_{\infty}/F_{\pr})\,.$$
By class field theory, the image of $A_{\pr}^*$ in $Gal(K_{\infty}/F_{\pr})$ is exactly $Gal(K_{\infty}/K)$. Let $a\in A$ be a generator of a prime ideal $\frak q\neq\pr$ and assume $sgn(a)=1$, so that $\Phi_a=\Phi_{\frak q}$; moreover, let $Frob_{\frak q}\in Gal(H^+(\Phi[\pr^\infty])/F)$ denote the Frobenius at $\frak q$. Passing from local to global class field theory one finds $Frob_{\frak q}=(a^{-1},K_{\infty}/F_{\pr})$. By \cite[Proposition 7.5.4]{goss} $\Phi_a$ acts on $\omega$ as $Frob_{\frak q}$ and hence $\chi^{-1}(a)=(a^{-1},K_{\infty}/F_{\pr})$.\\
By Tchebotarev for any $n$ all elements in $Gal(H^+(\Phi[\pr^n])/H^+)$ can be represented as $Frob_{\frak q}$ for some $\frak q$ as above. Therefore the corresponding $a$'s are dense in $A_{\pr}^*$ and we get
\begin{equation} \label{lcft-eq} (u,K_{\infty}/K)(\omega)=\Phi_{u^{-1}}(\omega) \end{equation}
for all $u\in A_{\pr}^*$.

\section{Coleman's formalism}\label{coleman}

Let $\mathcal K$ be a local field and $\{\mathcal K_n\}$ the tower of extensions of $\mathcal K$ generated by the torsion of a Lubin-Tate formal module: in \cite{col}, Coleman discovered an isomorphism between $\liminv\mathcal K_n^{\ast}$ and the subgroup of $\ol_{\mathcal K}((x))^{\ast}$ fixed under a certain operator $\n$. The same formalism can be used in the context of Drinfeld modules, as follows.

\subsection{A bit of functional analysis} Let $R$ be a subring of $\cp$: then, as usual, $R((x)):=R[[x]](x^{-1})$ is the ring of formal Laurent series with coefficients in $R$. Moreover, following \cite{col} we define $R[[x]]_1$ and $R((x))_1$ as the subrings consisting of those (Laurent) power series which converge on the punctured open ball
$$B':=B(0,1)-\{0\}\subset\cp\,.$$

These latter rings are endowed with a structure of topological $R$-algebras, induced by the family of seminorms $\{\|\cdot\|_r\}$, where $r$ varies in $|\cp|\cap(0,1)$ and $\|f\|_r:=\sup\{|f(z)|:|z|=r\}.$ One easily checks that this is the same as the ``compact-open" topology  of \cite[pag. 93]{col}: in particular a sequence $\{f_n\}$ in $\ol((x))_1$ converges to $f\in \ol((x))_1$ if and only if for each closed annulus $C$ around zero in $B'$ and for each $\epsilon>0$ there exists a positive integer $N(C,\epsilon)$ such that $|f_n(a)-f(a)|<\epsilon$ for all $a\in C$ and all $n\geq N(C,\epsilon)$. \\

Let $\phi\in \ol\{\tau\}$ be an additive polynomial having all its zeroes in $\mm$. One checks easily that, since $|\phi(z)|\leq |z|$ for all $z\in B'$, the map $\circ\phi: g\mapsto g\circ\phi$ defines a continuous endomorphism of $K[[x]]_1$.

\begin{lem} \label{lemcol} Let $\phi$ be as above; furthermore assume that $\phi(x)$ is separable. Then the image of $\circ\phi:K[[x]]_1\rightarrow K[[x]]_1$ consists exactly of those $g$ such that $g(x+u)=g(x)$ for all $u$ zeroes of $\phi$.
\end{lem}

This is essentially Lemma 3 of \cite{col}.

\begin{prf} For $u\in\mm$ let $T_u$ be the automorphism of $\cp[[x]]_1$ given by $g\mapsto g(x+u)$. The inclusion
$$\circ\phi(K[[x]]_1)\subset\bigcap_{\phi(u)=0}\ker(T_u-\id)$$ is clear.

Vice versa, assume that $f\in K[[x]]_1$ is $T_u$-invariant for all $u$'s. Let $K[[x]]_{\overline r}$ consist of those power series converging on the closed ball $\overline B(0,r)$: then $K[[x]]_1$ is the inverse limit of $K[[x]]_{\overline r}$ for $r<1$. The Weierstrass division Theorem holds in each $K[[x]]_{\overline r}$ and reasoning as in \cite[Lemma 3]{col} one can find $f_i$'s such that
$$f(x)=\sum_{i=0}^{n-1}f_i(0)\phi(x)^i+f_n(x)\phi(x)^n$$
for all $n\geq 0$. It remains to show that $f_n\phi^n$ tends to zero: one proves inductively that $\|f_n\|_r\leq\|f\|_r\,\|\phi\|_r^{-n}$, which implies, for $s>r$,
$$\|f_n\phi^n\|_r\leq\|f_n\|_s\,\|\phi^n\|_r\leq\|f\|_s\left(\frac{\|\phi\|_r}{\|\phi\|_s}\right)^n$$
because $\|g\|_r\leq\|g\|_s$\,. To conclude notice that $\|\phi\|_r<\|\phi\|_s$\,.
\end{prf}

\noindent {\bf Remark.} For the goals of this paper, it would have been enough to prove the weaker statement that $(\circ\phi)(K[[x]]_1)$ is the closure of the subspace of $T_u$-invariant polynomials. This can be done without Weierstrass theory, as follows.

Suppose $f\in K[x]$ is $T_u$-invariant for all the zeroes of $\phi$: we can assume inductively (taking as first step the constants) that if $g\in K[x]$ enjoys this property and $\deg(g)<\deg(f)$ then $g$ belongs to the image of $\circ\phi$. By the euclidean algorithm for $K[x]$, $f=f_1\phi+r$; evaluation in the zeroes of $\phi$ shows that $r=f(0)$ is a constant and then it is immediate to check that $f_1$ is $T_u$-invariant. Finally observe that $K[x]$ is dense in $K[[x]]_1$ and the maps are all continuous.

\begin{cor} \label{opmth} The map $\circ\phi$ induces an isomorphism of topological algebras between $K[[x]]_1$ and its image.
\end{cor}

\begin{prf} Observe that $K[[x]]_1$ is a Fr\'echet space over $K$ (for definitions and basic properties, see e.g. \cite[I,\S8]{schn}). Lemma \ref{lemcol} implies that $\circ\phi(K[[x]]_1)$, being closed, is Fr\'echet as well. Since $\circ\phi$ is injective, the corollary follows from the open mapping Theorem, as in \cite[Corollary
8.7]{schn}.
\end{prf}

\subsubsection{Some topological rings} \label{topologicalrings} To enhance clarity, we add a brief digression on topological structures for rings of power series. As above $R$ is a subring of $\cp$.

We let $R[[x]]\simeq R^{\N}$ be a topological $R$-algebra with the product topology: that is, a fundamental system of neighbourhoods of $0$ is given by
$$U_{\varepsilon,n}:=\left\{\sum a_ix^i\in R[[x]]:|a_i|<\varepsilon\;\forall\,i<n\right\}.$$
When we just write $R((x))$, we think of it as the additive group with the topology induced by the one on $R[[x]]$ (observe however that with this topology $R((x))$ is not exactly a topological ring, since multiplication by $x^{-1}$ is not continuous).

To compare structures remember that if $\{f_n\}$ converges to $f$ in $R[[x]]_1$ then the individual coefficients of the power series $f_n$ converge to those of $f$: it follows that the inclusion $R[[x]]_1\hookrightarrow R[[x]]$ is continuous. Observe, however, that $R((x))_1$ is not continuously injected in $R((x))\,$: e.g., if $|a|<1$ then $a^{n!}x^{-n}$ converges to $0$ in $\ol((x))_1$, but not in $\ol((x))\,$. This example shows as well that $R((x))_1$ is not complete.

When $R$ is a subring of the closed ball $\overline B(0,1)\subset\cp$, $R[[x]]_1\simeq R[[x]]$ as topological spaces and $R((x))_1=R((x))$ as sets. Moreover in this case we can furnish $R((x))$ with a third topology, defined (if the restriction of $v$ to $R$ is discrete) by the valuation $\nu(\sum a_ix^i):=\min_i\{v(a_i)\}.$ Once again, continuity of the inclusion $(R[[x]],\nu)\hookrightarrow R[[x]]_1$ fails to extend to Laurent series.

Lastly we remark that, for $\phi$ as above, $\circ\phi$ is a continuous endomorphism of $K((x))$ (but, of course, not of $K((x))_1$, unless $\phi$ has no zeroes in $B'$) and even an automorphism if $\phi(x)$ is separable. In fact $\phi(x)$ separable means that its degree 1 coefficient is not zero, hence one can find $\psi\in K[[x]]$ such that $\psi(\phi(x))=x$.

\subsection{Coleman's Theorems}

\begin{thm}[Coleman] \label{teo4col} There exist unique continuous operators
$$\T,\n:K((x))_1\rightarrow K((x))_1$$
such that respectively
$$\sum_{u\in\Phi[\pr^f]}g(x+u)=(\T g)\circ\Phi_{\eta}$$
$$\prod_{u\in\Phi[\pr^f]}g(x+u)=(\n g)\circ\Phi_{\eta}.$$
Moreover, $\T$ is a homomorphism of the additive group $K((x))_1$ and $\n$ of $K((x))_1^{\ast}$.  \end{thm}

\begin{prf} On $K[[x]]_1$ the theorem is an immediate consequence of Corollary \ref{opmth}: $\T$ (respectively $\n$) is just the composition of $(\circ\Phi_{\eta})^{-1}$ with $\sum T_u$ (resp.\! $\prod T_u$).\\
In order to extend $\T$ and $\n$ to $K((x))_1$, remember that $\Phi_{\eta}$ belongs to $x\ol[x]$: then, if $g\in K((x))_1$, for some $i\geq 0$ one has $\Phi_{\eta}(x)^ig\in K[[x]]_1$ and we put $\T(g):=x^{-i}\T(\Phi_{\eta}(x)^ig)$, $\n(g):=x^{-q_{\pr}^{if}}\n(\Phi_{\eta}(x)^ig)$. These are well-defined: e.g., if $g\in K[[x]]_1$,
$$\n(\Phi_{\eta}(x)g)\circ\Phi_{\eta}=\prod T_u(\Phi_{\eta}(x)g)=\prod T_u(\Phi_{\eta})T_u(g)=\Phi_{\eta}^{q_{\pr}^f}\prod T_u(g)=(x^{q_{\pr}^f}\n g)\circ\Phi_{\eta}.$$

\noindent Additivity of $\T$ and multiplicativity of $\n$ are immediate.
\end{prf}

As usual, we call $\T$ and $\n$ respectively the Coleman trace and norm.

\begin{lem} \label{nk} The equality $\displaystyle \n^kg\circ\Phi_{\eta}^k=\prod_{u\in\Phi[\pr^{fk}]}g(x+u)$ holds for any $g\in K((x))_1^{\ast}$.
\end{lem}

\begin{prf} Assume by induction that the statement is true up to $k-1$. Let $W\subset\Phi[\pr^{fk}]$ be a set such that $\Phi_{\eta}:W\rightarrow\Phi[\pr^{fk-f}]$ is a bijection. We have the following equalities:
$$\n^kg(\Phi_{\eta}^k(x))=\n^{k-1}\n g(\Phi_{\eta}^{k-1}(\Phi_{\eta}(x)))=$$
$$=\prod_{v\in\Phi[\pr^{fk-f}]}\n g(\Phi_{\eta}(x)+v)=\prod_{w\in W}\n g(\Phi_{\eta}(x+w))=$$
$$=\prod_{w\in W} (\n g\circ\Phi_{\eta})(x+w)=\prod_{W}\prod_{u\in \Phi[\pr^f]}g(x+u+w)=\prod_{t\in\Phi[\pr^{fk}]}g(x+t)\,.$$
\end{prf}

\begin{lem} \label{norm} One computes: $$(\mathcal N^kg)(\omega_n)=N_n^{n+k}(g(\omega_{n+k})).$$ Similarly,
$\T^kg(\omega_n)=Tr^{n+k}_n(g(\omega_{n+k}))$.
\end{lem}

\begin{prf} Replace $\omega_n=\Phi_{\eta}^k(\omega_{n+k})$ and apply Lemma \ref{nk}.
\end{prf}

\begin{lem} \label{ninfty} The restriction to $\ol((x))_1^{\ast}$ of the sequence of operators $\n^k$ converges to a continuous
endomorphism $\n^{\infty}$. \end{lem}

\begin{prf} Observe that $\ol((x))_1^{\ast}=x^{\Z}\times\ol[[x]]^{\ast}$ and that $\n x=x$. Therefore the lemma is proven if we show that, for any $g\in\ol[[x]]^{\ast}$, $\n^kg$ is a Cauchy sequence with respect to the valuation topology on $\ol[[x]]$ (uniformly on $g$). More precisely, we are going to prove that $\n^{k+1}g\equiv \n^kg\mod\m_K^{k+1}$ by induction on $k$.

First notice that
$$\Phi_{\eta}(x)=\prod_{u\in\Phi[\pr^f]}(x-u)\equiv x^{q_{\pr}^f}\mod {\frak m}_K$$
since $v(u)>0$ for $u\in\Phi[\pr^f]$. Therefore
$$g^{q_{\pr}^f}\equiv\prod g(x+u)=\n g\circ\Phi_{\eta}\equiv\n g(x^{q_{\pr}^f})\equiv(\n g(x))^{q_{\pr}^f}\mod {\frak m}_K$$
where the last congruence is true because $q_{\pr}^f=|\ol_K/\m_K|$. It follows that $\n g\equiv g\mod {\frak m}_K$.

Now put $h:=\frac{\n^kg}{\n^{k-1}g}\,$, so that our claim becomes $\n h\equiv 1\mod\m_K^{k+1}$. By the induction hypothesis $h=1+\pi^kg_1$ and this implies
$$(\n h\circ\Phi_{\eta})(x)=\prod_{u\in\Phi[\pr^f]}(1+\pi^kg_1(x+u))\equiv (1+\pi^kg_1(x))^{q_{\pr}^f}\equiv 1\mod{\frak m}^{k+1}_K.$$
To conclude, observe that since $\Phi_{\eta}$ is monic $\nu(h\circ\Phi_{\eta})=\nu(h)$ for any $h$ (where $\nu$ is the valuation on $\ol[[x]]$ defined in \S\ref{topologicalrings}).
\end{prf}

Of course $\n\circ\n^{\infty}=\n^{\infty}$ and $\n^{\infty}$ is a projection.

\begin{thm}[Coleman] \label{teo8} The evaluation map $ev:f\mapsto\{f(\omega_n)\}$ gives an isomorphism
$$(\ol((x))^{\ast})^{\n=id}\simeq\liminv K_n^{\ast}$$
where the inverse limit is taken with respect to the norm maps.
\end{thm}

\begin{prf} The map is injective, because a function is uniquely determined by its values at the $\omega_n$'s (e.g., observe that $|\omega_n|<1$ and use \cite[Proposition 2.11]{goss}).

Notice that $\ol((x))^{\ast}=x^{\Z}\times\ol[[x]]^{\ast}$ and $\liminv K_n^{\ast}=\omega^{\Z}\times\liminv\ol_n^{\ast}$: since $ev(x)=\omega$, it suffices to show $(\ol[[x]]^{\ast})^{\n=\id}\simeq\liminv\ol_n^{\ast}$.\\
Consider the diagram
\[ \begin{CD}
\ol[[x]]^{\ast} @>ev>> \prod\ol_n^{\ast}\\
  @VV{\n/\id}V     @VV{N/\id}V\\
\ol[[x]]^{\ast} @>ev>> \prod \ol_n^{\ast}
 \end{CD} \]\\
where $N$ is the norm map $(x_n)\mapsto(N_n^{n+1}x_{n+1})$. It commutes by Lemma \ref{norm};
$\liminv\ol_n^{\ast}$ is the kernel of the right-hand side:
hence $ev(g)\in\liminv\ol_n^{\ast}$ iff $\n g=g$. Since $(\ol[[x]]^{\ast})^{\n=\id}=\n^{\infty}(\ol[[x]]^{\ast})$
is compact (because so is $\ol[[x]]^{\ast}$),
the theorem is proven if we show that the image of $ev$ is
dense in a set containing $\liminv\ol_n^{\ast}$.
For any $u=(u_n)\in\liminv\ol_n^{\ast}$ and any $k$ there exists $g\in\ol[[x]]^*$ such that $g(\omega_{2k})=u_{2k}$. Let $h:=\n^kg$. Remembering (from the proof of Lemma \ref{ninfty}) that $\n^kg\equiv\n^{k+r}g$ mod $\m^k$ for any $r\geq 0$, we get
$$h(\omega_i)\equiv\n^{2k-i}g(\omega_i)=N^{2k}_i(g(\omega_{2k}))=u_i\mod\m^k$$
for all $i=1,...,k$: density follows.
\end{prf}

In particular, $\omega\in\liminv\ol_n^{\ast}$ corresponds to $x\in\ol((x))^{\ast}$ and, more generally, for $a\in A_{\pr}^*$ the element $a\cdot\omega$ corresponds to $\Phi_a(x)\in\ol((x))^*$ (this is equivalent to changing generator of $T_{\pr}\Phi$).\\

We conclude with a lemma we are going to use in the next section.

\begin{lem} \label{dlog} Let $\dlog:K((x))_1^{\ast}\rightarrow K((x))_1$ be the logarithmic derivative operator, $\dlog(g):=\frac{g'}{g}$. Then
$$\T\dlog\, g=\eta\dlog\n g.$$
\end{lem}

\begin{prf} One just computes:
$$(\T\dlog g)\circ\Phi_{\eta}=\sum T_u(\dlog g)=\dlog\prod T_u(g)=\dlog(\n g\circ\Phi_{\eta})=\eta(\dlog(\n g)\circ\Phi_{\eta})$$
using the fact that $\dlog$ is a homomorphism and $\frac{d}{dx}\Phi_{\eta}(x)=\eta$.
\end{prf}

\subsection{Higher rank Drinfeld modules} As the reader may have noticed, the rank of the Drinfeld module $\Phi$ plays essentially no part in Theorem \ref{teo4col}. This suggests the possibility of extending Coleman's results to Drinfeld modules of any rank: we sketch an approach.

In this subsection (and only here) our notations are slightly modified. For simplicity we take $F:=\F_q(T)$ and $A:=\F_q[T]\,$: then, fixing a prime ideal $\pr=(\pi)$ of $A$, we have $K=F_{\pr}$ and $\ol=A_{\pr}$. Let $\Phi:A\rightarrow\ol\{\tau\}$ be a rank $r$ Drinfeld module, $r>1$; we need the hypothesis that $\Phi$ has good reduction mod $\m_K$ and that the height is maximal: $h=r$. In particular, it follows that all zeroes of $\Phi_{\pi}$ are in the open unit ball $B(0,1)\subset{\bf C}_{\pr}$. We also assume that $\Phi_{\pi}(x)$ is a monic polynomial. Then, reasoning exactly as in Theorem \ref{teo4col} and Lemma \ref{ninfty}, one proves the following.

\begin{thm} There exists a continuous homomorphism $\n: K((x))_1^*\rightarrow K((x))^*_1$ such that
$$\prod_{u\in\Phi[\pr]}g(x+u)=(\n g)\circ\Phi_{\pi}.$$
The restriction to $\ol((x))_1^{\ast}$ of the sequence of operators $\n^k$ converges to a continuous endomorphism $\n^{\infty}$.
\end{thm}

As in section \ref{tower}, we choose a sequence $\{\omega_n\}_{n\geq 1}$ so that
$$\Phi_{\pi}^n(\omega_n)=0\neq\Phi_{\pi}^{n-1}(\omega_n)\text{ and }\Phi_{\pi}(\omega_{n+1})=\omega_n$$
and construct a tower $\{K_n\}$ by $K_n:=K_{n-1}(\omega_n)$, with $K_0:=K$. Because of the rank, these extensions are much smaller than $K(\Phi[\pr^n])$ and they are not Galois; however, since the polynomials $\Phi_{\pi}(x)x^{-1}$ and $\Phi_{\pi}(x)-\omega_n$ are Eisenstein, it still holds that each $K_n/K$ is totally ramified, with uniformizer $\omega_n$.\\
For any $n\geq 1$ there is a norm map
$$N^{n+1}_n:K_{n+1}^*\rightarrow K_n^*,\hspace{7pt}a\mapsto\prod_{\sigma \in S_{n+1}}\sigma(a)$$
where the product is taken on the set of embeddings
$$S_{n+1}:=\{\sigma:K_{n+1}\hookrightarrow{\bf C}_{\pr}:\sigma|_{K_n}=id_{K_n}\}\,.$$
It follows from the additivity of $\Phi_{\pi}$ that the assignment $\sigma\mapsto\sigma(\omega_{n+1})-\omega_{n+1}$ is a bijection $S_{n+1}\rightarrow\Phi[\pr]\,$: therefore $N^{n+1}_n(\omega_{n+1})=\omega_n$ and, more generally,
$$(\n g)(\omega_n)=\n g(\Phi_{\pi}(\omega_{n+1}))=\prod_{u\in\Phi[\pr]}g(\omega_{n+1}+u)=N^{n+1}_n(g(\omega_{n+1}))\,.$$

\begin{thm} The evaluation map $ev:f\mapsto\{f(\omega_n)\}$ gives an isomorphism
$$(\ol((x))^{\ast})^{\n=id}\simeq\liminv K_n^{\ast}\,.$$
\end{thm}

\noindent The proof is the same as for Theorem \ref{teo8}.

\section{The explicit reciprocity law} \label{law}

The reader is reminded that $Tr_n$, $N_n$ denote respectively trace and norm from $K_n$ to $F_{\pr}$. Also, we let
$$(\cdot, L^{ab}/L):L^{\ast}\longrightarrow G_{L}^{ab}$$
be the local norm symbol map and write $Col_u$ for the power series in $\ol((x))^{\ast}$ associated to $u\in\liminv K_n^{\ast}$ by Coleman's isomorphism of Theorem \ref{teo8}. To lighten notation, in this section the action of $A_{\pr}$ via $\Phi$ will be often denoted by $a\cdot x:=\Phi_a(x)$.

\subsection{The Kummer pairing} Next to $\liminv K_n^*$ we consider $\limdir K_n$, defined as the direct limit of the maps $\Phi_{\eta}: K_n\rightarrow K_{n+1}$: that is, $\limdir K_n$ consists of sequences $a=(a_n)_{n\geq N}$ (for some $N\in\N$) such that $a_N\in K_N$ and $a_{n+1}=\Phi_{\eta}(a_n)$, modulo the relation $(a_n)_{n\geq N}=(b_n)_{n\geq M}$ if $a_n=b_n$ for $n\gg 0$.

Kummer theory yields a pairing $(\;,\;)_n:K_n\times K_n^{\ast}\rightarrow\Phi[\pr^{fn}]$ defined by
$$(a,u)_n:=((u,K_n^{ab}/K_n)-1)(\,\psqrt[n]{a}).$$
Here $\psqrt[n]{a}$ is a solution of $\Phi_{\eta}^n(x)=a$: since any two roots differ by an element in $\Phi[\pr^{fn}]\subset K_n$, the value of $(a,u)_n$ is independent of the choice of $\psqrt[n]{a}$.

Observe that, since by definition $\psqrt[n]{a}={}_{\eta}\!\!\!\!\!\!\!\sqrt[n+1]{\Phi_{\eta}(a)}$,
$$(\Phi_{\eta}(a),u)_{n+1}=((u,K_{n+1}^{ab}/K_{n+1})-1)(\;{}_{\eta}\!\!\!\!\!\!\!\sqrt[n+1]{\Phi_{\eta}(a)})=$$
$$=((N^{n+1}_n(u),K^{ab}_n/K_n)-1)(\,\psqrt[n]{a})=(a,N^{n+1}_n(u))_n \,.$$
This means that, given $a=(a_n)\in\limdir K_n$ and $u=(u_n)\in\liminv K_n^*$, one has (for any $n$ large enough that $a_n$ exists)
$$(a_{n+1},u_{n+1})_{n+1}=(a_n,u_n)_n\,.$$
Therefore we can define a limit form of the Kummer pairing
$$(\;,\;):\limdir K_n\times\liminv K_n^{\ast}\rightarrow\Phi[\pr^{\infty}]$$
by $(a,u):=(a_n,u_n)_n$ for $n\gg0$.\\

One checks immediately that $(\;,\;)$ is bilinear, additive in the first variable and multiplicative in the second. In particular, since values are in a group of exponent $p$, it follows that $(\cdot,\zeta)_n= 0$ for any root of unity $\zeta\in K_n^*\,.$\\

\begin{lem} \label{cont} All the pairings $(\;,\;)_n$ are continuous.
Furthermore, $(a,\cdot)_n\equiv 0$ for any $a\in K_n$ such that $v(a)>nf+(q_\pr-1)^{-1}$.
\end{lem}

\begin{prf} For any $a\in K_n$, $(a,\cdot)_n$ is continuous: therefore the first assertion follows from the second, which in turn is an easy application of Krasner's Lemma. In fact, if one can choose $\alpha=\psqrt[n]{a}$ so that $v(\alpha)$ is big enough, then it follows $(a,u)=0$ because $|(a,u)|\leq |\alpha|=|(u,K_n^{ab}/K_n)(\alpha)|$ and $\Phi[\pr^{nf}]$ is a discrete subset of $\dc$.

We are left with a valuation computation. First of all, observe that for any $j\geq 1$ if $u\in\Phi[\pr^j]-\Phi[\pr^{j-1}]$ then $$v_j:=v(u)=[K(\Phi[\pr^j]):K]^{-1}=|(A/\pr^j)^*|^{-1}=\frac{1}{q_\pr^{j-1}(q_\pr-1)}$$
(because $\frac{\Phi_{\pr^j}(x)}{\Phi_{\pr^{j-1}}(x)}$ is Eisenstein). In particular the smallest non-zero elements in $\Phi[\pr^{nf}]$ have valuation $v_1=(q_\pr-1)^{-1}$. We also put $v_0:=\infty$.

Now choose $\alpha$ a root of $\Phi_{\eta}^n(x)=a$ such that $v(\alpha)$ is maximal: we get $$a=\prod_{u\in\Phi[\pr^{fn}]}(\alpha+u)\; \text{ and }\; v(a)=\sum_{u\in\Phi[\pr^{fn}]}v(\alpha+u)\,,$$
with $v(\alpha+u)=\min\{v(\alpha),v(u)\}$ for all $u$'s because of the maximality hypothesis. Hence if $v_j\geq v(\alpha)>v_{j+1}$ we obtain
$$v(a)=\sum_{u\in\Phi[\pr^j]}v(\alpha)+\sum_{u\in\Phi[\pr^{nf}]-\Phi[\pr^j]}v(u)=q_{\pr}^jv(\alpha)+nf-j$$
and $1+nf-j+(q_\pr-1)^{-1}=q_{\pr}^jv_j+nf-j\geq v(a)>nf-j+(q_\pr-1)^{-1}$.
\end{prf}

For a more detailed analysis of how $v(a)$ determines the extension $K_1(\psqrt[n]{a})/K_1$ see \cite[Proposition 2.1]{angles}.\\
We remark that the computation in the proof of Lemma \ref{cont} yields immediately the following result.

\begin{lem} \label{nf-c} For any $a=(a_n)_{n\geq N}\in\limdir\m_n$ there exists a constant $c(a)$ such that $v(a_n)\geq nf-c(a)$ for all $n\geq N$.
\end{lem}

\begin{lem} \label{cl1} Let $a\in K_n^*$: then $(a,a)_n=0$.
\end{lem}

\begin{prf} Let $\alpha$ be a representative of $\psqrt[n]{a}$ and put $L:=K_n(\alpha)$. Kummer theory identifies $\Gal(L/K_n)$ with a subgroup $V$ of $\Phi[\pr^{nf}]$: then one sees that
$$a=\prod_{u\in \Phi[\pr^{nf}]/V}\prod_{v\in V}(\alpha+u+v)=\prod_{\Phi[\pr^{nf}]/V}N_{L/K_n}(\alpha+u)$$
and consequently $(a,K^{ab}_n/K_n)$ acts trivially on $L$.
\end{prf}

\begin{lem} \label{1-x} Let $b\in\m_n-\{0\}$: then $(c,1-b)_n=(\frac{bc}{1-b},b^{-1})_n$ for all $c\in K_n$.
\end{lem}

\begin{prf} If $c=0$ both sides are 0. If not, by applying Lemma \ref{cl1} to $a=c(1-b)$ we get, by bilinearity,
$$(c,1-b)_n=(cb,c)_n+(cb,1-b)_n$$
and, by recurrence,
$$(c,1-b)_n=\sum_{j=1}^{\infty}(cb^j,cb^{j-1})_n$$
(the sum converges because only a finite number of terms are not 0, by Lemma \ref{cont}).\\
By Lemma \ref{cl1} we have
$(cb^j,cb^{j-1})_n=(cb^j,cb^j)_n+(cb^j,b^{-1})_n=(cb^j,b^{-1})_n$
and therefore
$$(c,1-b)_n=\sum_{j=1}^{\infty}(cb^j,b^{-1})_n=(c\sum_{j=1}^{\infty}b^j,b^{-1})_n=(\frac{cb}{1-b},b^{-1})_n\,.$$
\end{prf}

\begin{prop} \label{kumdlog} Let $a\in\limdir\m_n$, $u\in\liminv K_n^*$: then
$$(a,u)=(a_n\omega_n\dlog Col_u(\omega_n),\omega_n)_n$$
for all $n$ sufficiently large.
\end{prop}

\begin{prf} The statement is clearly true for $u=\omega$. As for $u\in\liminv\ol_n^*$, we are going to prove that, more generally,
$$(c,w)_n=\big(c\omega_n\frac{\dlog w}{d\omega_n},\,\omega_n\big)_n$$
for all $(c,w)\in \m_n^{s_n}\times\ol_n^*$ for some $s_n$. Here
$\dlog:\ol_n^*\rightarrow\Omega_{\ol_n/\ol}$ is the map
$x\mapsto\frac{dx}{x}$. The module of differentials is free over
$\ol_n/\frak D_n$ with generator $d\omega_n$: by Lemma \ref{lema3}
if $c\in\m_n$ with $v(c)>\frac{2}{q_{\pr}-1}-v(\omega_n)$ and
$\delta\in\frak D_n$ the inequality
$$v(c\omega_n\delta)\geq v(c)+v(\omega_n)+nf-\frac{1}{q_{\pr}-1}>nf+\frac{1}{q_{\pr}-1}$$
holds and so Lemma \ref{cont} shows that $(c\omega_n\frac{\dlog
w}{d\omega_n},\omega_n)_n$ is well-defined. Thanks to Lemma
\ref{nf-c} this is enough for our purposes.

Observe that it suffices to prove the claim for $w=1-\zeta\omega_n^k$ ($\zeta$ varying among roots of unity in $K_n$), because one can choose a topological basis of $1+\m_n$ consisting of elements of this form and the pairing is continuous and linear.\\
Applying Lemma \ref{1-x} with $x=\zeta\omega_n^k$ we get
$$(c,1-\zeta\omega_n^k)_n=(\frac{c\zeta\omega_n^k}{1-\zeta\omega_n^k},\omega_n^{-k})_n=(c\omega_n\frac{\dlog(1-\zeta\omega_n^k)}{d\omega_n},\omega_n)_n$$
because
$$\dlog(1-\zeta\omega_n^k)=\frac{-k\zeta\omega_n^{k-1}}{1-\zeta\omega_n^k}d\omega_n\,.$$
To conclude just notice that, for $u\in\liminv K_n^*$, $\displaystyle\dlog Col_u(\omega_n)=\frac{\dlog u_n}{d\omega_n}\,.$ (Caveat! In this last formula the symbol $\dlog$ appears with two different meanings: on the left-hand side $\dlog Col_u$ is the power series $\frac{1}{Col_u(x)}\frac{d}{dx}Col_u(x)$, evaluated in $\omega_n$, while on the right-hand side $\dlog:\ol_n^*\rightarrow\Omega_{\ol_n/\ol}$ is the map we defined above.)
\end{prf}

\subsection{The analytic pairing} As above, let $u\in\liminv K_n^*$. Lemmata \ref{norm} and \ref{dlog} together with the $\n$-invariance of $Col_u$ yield
$$Tr^{n+k}_n\dlog Col_u(\omega_{n+k})=\eta^k\dlog Col_u(\omega_n).$$
Besides, for $a=(a_n)\in\limdir\frak m_n$,
$\lambda(a_{n+k})=\eta^k\lambda(a_n)$ by definition of $\lambda$
(Proposition \ref{lambda}). It follows that
$$Tr_{n+k}(\eta^{-n-k}\lambda(a_{n+k})\dlog Col_u(\omega_{n+k}))=\eta^kTr_n(\eta^{-n}\lambda(a_n)\dlog Col_u(\omega_n)).$$

\begin{lem} \label{2pair} $Tr_n(\eta^{-n}\lambda(a_n)\dlog Col_u(\omega_n))\in A_{\pr}$ for $n\gg 0$. \end{lem}

\begin{prf} By Lemma \ref{nf-c}, $v(a_n)\geq nf-c$ for some constant $c$ (depending on $a$): Proposition \ref{lambda} implies that the same is true for $\lambda(a_n)$. Since $Col_u\in x^{\Z}\times\ol[[x]]^*$ one has $v(\dlog Col_u(\omega_n))\geq -v(\omega_n)$ and the last value, in turn, is controlled by Lemma \ref{lema3}. Now apply Corollary \ref{trace}.
\end{prf}

\noindent Therefore one can define a second pairing
$$[\;,\;]:\limdir \frak m_n\times\liminv K_n^{\ast}\rightarrow \Phi[\pr^{\infty}]$$
putting
$$[a,u]:=Tr_n(\eta^{-n}\lambda(a_n)\dlog Col_u(\omega_n))\cdot\omega_n$$
for $n\gg 0$. (Recall that, by definition, $a\cdot\omega_n=\Phi_a(\omega_n)\,$.)

\subsubsection{} \label{anpairn} It is convenient to define also a level $n$ pairing $[\;,\;]_n:\frak m_n^{t_n}\times K_n^*\rightarrow \Phi[\pr^{nf}]\,$ for some $t_n\geq 1$.\\
The logarithmic differential we used in the proof of Proposition
\ref{kumdlog} can be extended to a homomorphism
$$\frac{\dlog}{d\omega_n}:K_n^*=\omega_n^{\Z}\times\ol_n^*\rightarrow\m_n^{-1}/\frak D_n$$
by putting $\displaystyle\frac{\dlog\omega_n^i}{d\omega_n}:=i\omega_n^{-1}\,.$

\begin{lem}\label{lema20} The pairing
$$[a,u]_n:=Tr_n\left(\eta^{-n}\lambda(a)\frac{\dlog u}{d\omega_n}\right)\cdot\omega_n$$
is well defined for $v(a)\geq 2(q-1)^{-1}$.
\end{lem}

\begin{prf} We have to show that $Tr_n(\eta^{-n}\lambda(a)b)$ belongs to $A_{\pr}$ for any $b\in\m_n^{-1}$ and that
$$v\left(Tr_n(\eta^{-n}\lambda(a)\delta)\right)\geq nf$$
for $\delta\in\frak D_n\,.$ Since the hypothesis implies
$v(\lambda(a))=v(a)$, both assertions are easy consequences of
Lemma \ref{lema3} and Corollary \ref{trace}.
\end{prf}

It is clear that if $(a,u)\in\limdir\m_n\times\liminv K_n^*$, the equality $[a,u]=[a_n,u_n]_n$ holds for $n\gg 0$.

\begin{prop} \label{andlog} For $n$ large enough, $[a,u]=[a_n\omega_n\dlog Col_u(\omega_n),\omega_n]_n$.
\end{prop}

\begin{prf} To lighten notation, put $D=\omega_n\dlog Col_u(\omega_n)$; observe that $D\in\ol_n$. One has to check that
$$Tr_n(\eta^{-n}\lambda(a_n)D\omega_n^{-1})\cdot\omega_n=Tr_n(\eta^{-n}\lambda(a_nD)\omega_n^{-1})\cdot\omega_n,$$
i.e. that
$$v\big(Tr_n\left(\frac{\lambda(a_nD)-\lambda(a_n)D}{\omega_n}\right)\big)\geq 2nf\,.$$
By Proposition \ref{lambda},
$$v\big(\lambda(a_nD)-\lambda(a_n)D\big)=v\big(\sum_{i=1}^{\infty} c_ia_n^{q^i}(D^{q^i}-D)\big)\geq\min_{i\geq 1}\{q^iv(a_n)-i\}$$
which for $n\gg 0$ is at least $2nf-c$ for some $c$ independent of $n$, by Lemma \ref{nf-c}. Now apply Corollary \ref{trace} to get rid of $c$.
\end{prf}

\subsection{The reciprocity law} \label{prooflaw} In this paragraph we prove Wiles' explicit reciprocity law for rank 1 Drinfeld modules (Theorem \ref{reclaw}).

\begin{prop} \label{receps} Let $a\in\frak m_n$, $v(a)\geq 2(q-1)^{-1}$: then $[a,\omega_n]_n=(a,\omega_n)_n\,.$
\end{prop}

\begin{prf} Our proof is divided in many steps, along the lines of \cite[I, \S4]{desh}.\\

\noindent 1. To start with, we put $a_m:=\Phi_{\eta}^{m-n}(a)$ for
all $m\geq n$ and let $b_m:=a_m\omega_m^{-1}$. Thanks to Lemmata
\ref{nf-c} and \ref{lema3} there is a constant $c$ (depending only
on $a$) such that $v(b_m)\geq mf-c$.

As a consequence of Lemma \ref{cl1} we get
$$0=(a_m+\omega_m,(1+b_m)\omega_m)_m=(a_m,\omega_m)_m+(a_m,1+b_m)_m+(\omega_m,1+b_m)_m$$
for $m\geq n$.\\

\noindent 2. {\em Claim:} if $m\gg 0$, then $(a_m,1+b_m)_m=0.$\\
{\em Proof.} By definition of the Kummer pairing, for $m>n$ we have
$$(a_m,1+b_m)_m=(a_n,N^m_n(1+b_m))_n.$$
Since $1+b_m$ tends to 1, so does also $N_n^m(1+b_m)\,$: the claim is proven because $(a,\cdot)_n$ is continuous and $\Phi[\pr^{nf}]$ discrete.\\

\noindent 3. {\em Claim:}
$(a_n,\omega_n)_n=\omega_{2m}-\Phi_{N_m(1+b_m)^{-1}}(\omega_{2m})$ for $m\gg 0$.\\
{\em Proof.} From the above, we have
$$(a_n,\omega_n)_n=(a_m,\omega_m)_m=-(\omega_m,1+b_m)_m=(1-(1+b_m,K_m^{ab}/K_m))(\,\psqrt[m]{\omega_m}).$$
We can take $\psqrt[m]{\omega_m}=\omega_{2m}$. The extension $K_{2m}/F_{\pr}$ is abelian: hence
$$(1+b_m,K_m^{ab}/K_m)_{|K_{2m}}=(N_m(1+b_m),F_{\pr}^{ab}/F_{\pr})_{|K_{2m}}\,.$$
Now one applies formula (\ref{lcft-eq}).\\

\noindent 4. {\em Claim:} $N_m(1+b_m)^{-1}\equiv 1-Tr_m(b_m)\mod\pr^{2mf}$ for $m$ big enough.\\
{\em Proof.} Take $k\in\N$ such that $kf>c+1$ where $c$ is the constant which appeared in step 1. We can assume $m\gg k$.\\
Put $\beta:=Tr^m_{m-k}(b_m)$: then $v(\beta)\geq
v(b_m)+kf-v(\omega_{m-k})$ by Corollary \ref{trace}. We have
$$ N^m_{m-k}\left(\frac{1}{1+b_m}\right)=1-\beta+\delta\,,$$
where $\delta$ is an element of $K_{m-k}$ such that $v(\delta)\geq 2v(b_m)$. Now apply $N_{m-k}$ to the above equation to obtain
$$ N_m(1+b_m)^{-1}=N_{m-k}(1-\beta+\delta)=1-Tr_{m-k}(\beta-\delta)+\theta\equiv 1-Tr_m(b_m)\mod\pr^{2mf},$$
because summands in $\theta$ have valuation at least
$$2v(\beta-\delta)\geq 2\min\{(m+k)f-c-v(\omega_{m-k}),2mf-2c\}$$
and $v(Tr_{m-k}(\delta))\geq\lfloor(3m-k)f-2c-(q_{\pr}-1)^{-1}\rfloor.$ \\

\noindent 5. From steps 3 and 4 we get $(a,\omega_n)_n=Tr_m(a_m\omega_m^{-1})\cdot\omega_{2m}\,.$ On the other hand $[a,\omega_n]_n=[(a_m),\omega]=[a_m,\omega_m]_m\,.$ Since $Col_{\omega}=x$, $\dlog Col_{\omega}=\frac{1}{x}$ and by
definition we obtain
$$[a_m,\omega_m]_m=Tr_m\!\!\left(\frac{\lambda(a_m)}{\eta^m\omega_m}\right)\cdot\omega_m=\frac{1}{\eta^m}Tr_m\!\!\left(\frac{\lambda(a_m)}{\omega_m}\right)\cdot(\eta^m\cdot\omega_{2m})=Tr_m\!\!\left(\frac{\lambda(a_m)}{\omega_m}\right)\cdot\omega_{2m}\,.$$
The proof is completed by the same reasoning as in Proposition \ref{andlog}.
\end{prf}

Combining Propositions \ref{kumdlog}, \ref{andlog} and \ref{receps}, we get our reciprocity law:

\begin{thm} \label{reclaw} The two pairings $(\;,\;)$ and $[\;,\;]$ on $\limdir \frak m_n\times\liminv K_n^{\ast}\rightarrow \Phi[\pr^{\infty}]$, defined respectively by
$$[a,u]:=Tr_n(\eta^{-n}\lambda(a_n)\dlog Col_u(\omega_n))\cdot\omega_n$$
and
$$(a,u):=((u_n,K_n^{ab}/K_n)-1)(\,\psqrt[n]{a_n})$$
for $n\gg0$, are equal. \end{thm}

\subsubsection{The Kummer pairing in practice}

A weaker form of our explicit reciprocity law can be used to calculate the Kummer pairing also when $u_n\in K_n^*$ is not a coordinate in an inverse limit or, even if it is the case, one does not know how to explicitly find $Col_u$.

Given $u_n\in K_n^*$ and $a_n\in\mathfrak{m}_n$ we want to compute
$(a_n,u_n)_n\,$. We need to impose that there exists $u_m$ for
some convenient $m\geq n$ such that $N^m_n(u_m)=u_n$. If so we
have:
$$(a_n,u_n)_n=(\Phi_{\eta}(a_n),u_{n+1})_{n+1}=\ldots=(\Phi_{\eta}^j(a_n),u_{n+j})_{n+j}$$
for any integer $0\leq j\leq m-n$. Put $a_{n+j}:=\Phi_{\eta}^j(a_n)$.

In particular suppose we can take $m$ big enough to have $v(a_m)>
2/(q-1)$ (and hence $v(a_m)>\frac{2}{q_{\pr}-1}-v(\omega_n)$). (An
estimate on the required size of $m-n$ can be obtained from the
computations proving Lemma \ref{nf-c}.) Then
$$(a_m,u_m)_m=(a_m\omega_m\frac{\dlog u_m}{d\omega_m},\omega_m)_m$$
(see the proof of Proposition \ref{kumdlog}) and by Proposition
\ref{receps} we have
$$(a_m,u_m)_m=(a_m\omega_m\frac{\dlog u_m}{d\omega_m},\omega_m)_m=[a_m\omega_m\frac{\dlog u_m}{d\omega_m},\omega_m]_m\,.$$
By Proposition \ref{lambda} we have
$v(\lambda(z)-z)\geq\min_{i\geq 1}\{q^iv(z)-i\}\,;$ in particular
this minimum is attained in $i=1$ when $v(z)>2(q-1)^{-1}$. Hence a
simple computation shows that the further condition
$v(a_m)\geq\frac{1}{q}\big(mf+1+(q_{\pr}-1)^{-1}+v(\omega_m)\big)$
implies
$$[a_m\omega_m\frac{\dlog u_m}{d\omega_m},\omega_m]_m=Tr_m\left(\frac{a_m}{\eta^m}\frac{\dlog u_m}{d\omega_m}\right)\cdot\omega_m=[a_m,u_m]_m\,.$$

The limit form of the Kummer pairing, as in Theorem \ref{reclaw},
is useful rather for the purposes of Iwasawa theory (which is not
yet well understood in our setting) than for a concrete
calculation of the pairing at a level $n$. Historically, the
interest in computing the Kummer pairing for local fields at a
finite level was originated by the study of diophantine equations,
in particular the Fermat one. For a survey on various explicit
reciprocity laws in the local case we refer to \cite{balo}.

\section*{Acknowledgements}

The first idea of this paper was born when the authors met in Universit\"at M\"unster, during fall 2000
(the first author was supported by SFB 478, the second one by the TMR network Arithmetic Algebraic Geometry
for mobility of young researchers).
The project was later carried on with support and/or hospitality by Dipartimento di Matematica of Universit\`a di Padova, Departament de Matem\`atiques in the Universitat Aut\`onoma de Barcelona, Dipartimento di Matematica of Universit\`a degli Studi di Milano and also Dipartimento di Matematica of Universit\`a di Pavia: we would like to thank all of these institutions. Thanks also to Matteo Longo, for his comments on a first draft of this paper. Finally we thank the referee for his or her comments and suggestions.

\vspace{1cm}

Francesc Bars Cortina, Depart. Matem\`atiques, Universitat Aut\`onoma de Barcelona, 08193 Bellaterra. Catalonia. Spain.\\
E-mail: francesc@mat.uab.cat \\

Ignazio Longhi,  Department of Mathematics, National Taiwan University, No. 1 section 4 Roosevelt Road, Taipei 106, Taiwan.\\
E-mail: longhi@math.ntu.edu.tw


\begin{thebibliography}{30}

\bibitem{angles} B. Angl\`es: On explicit reciprocity laws for the local Carlitz-Kummer symbols. {\em J. Number Theory} {\bf 78} (1999), no. 2, 228--252.
\bibitem{artinhasse} E. Artin, H. Hasse: Die beiden Erg\"anzungss\"atze zum Reziprozit\"atsgesetz der $l^n$-ten \-Po\-tenz\-re\-ste im K\"orper der $l^n$-ten Einheitswurzeln. {\em Abh. Math. Sem. Univ. Hamburg} {\bf 6} (1928), 146-162.
\bibitem{balo} F. Bars, I. Longhi: Reciprocity laws \`a la Iwasawa-Wiles. To appear in {\em Bibl. Rev. Mat. Iberoamericana}.
\bibitem{col} R. Coleman: Division values in local fields. {\em Invent. Math.} {\bf 53} (1979) 91-116.
\bibitem{desh} E. de Shalit:  Iwasawa theory of elliptic curves with complex multiplication. Perspectives in Mathematics, 3. Academic Press, Inc., Boston, MA, 1987.
\bibitem{goss} D. Goss: Basic structures of function field arithmetic. Springer-Verlag, New York, 1996.
\bibitem{hayes} D. Hayes: A brief introduction to Drinfeld modules. In: The arithmetic of function fields (Columbus, OH, 1991), 1-32. de Gruyter, Berlin, 1992.
\bibitem{iwasawa} K. Iwasawa: On explicit formulas for the norm residue symbol. {\em J. Math. Soc. Japan} {\bf 20} (1968), 151-164.
\bibitem{kato} K. Kato: Generalized explicit reciprocity laws. {\em Adv. Stud. Contemp. Math. (Pusan)} {\bf 1} (1999) 57-126
\bibitem{lang} S. Lang: Cyclotomic Fields I and II. GTM 121, Springer-Verlag, New York/Berlin, 1990.
\bibitem{rosen} M. Rosen:  Formal Drinfeld modules. {\em J. Number Theory} {\bf 103} (2003), no. 2, 234-256.
\bibitem{schn} P. Schneider: Nonarchimedean functional analysis. Springer-Verlag, Berlin, 2002.
\bibitem{scholl} A.J. Scholl: An introduction to Kato's Euler systems. In: Galois representations in arithmetic algebraic geometry (A.J. Scholl and R.L. Taylor eds.) {\em London Math. Soc. Lecture Notes} {\bf 254} (1998) 379-460
\bibitem{serre} J.-P. Serre: Local Fields. GTM 67, Springer-Verlag, New York/Berlin, 1979.
\bibitem{wash} L. Washington: Introduction to cyclotomic fields. Second edition. GTM 83, Springer-Verlag, New York, 1997
\bibitem{wiles} A. Wiles: Higher explicit reciprocity laws. {\em Ann. of Math.} {\bf 107} (1978), no. 2, 235-254.

\end{thebibliography}
\end{document}